# The Classical Smarandache Function and a Formula for Twin Primes

Dhananjay P. Mehendale
Sir Parashurambhau College, Tilak Road, Pune-411009, India.

## Abstract

This short paper presents an exact formula for counting twin prime pairs less than or equal to x in terms of the classical Smarandache Function. An extension of the formula to count prime pairs (p, p+2n), n > 1, is also given.

## 1. Introduction:

The most known Smarandache function which has become a classical Smarandache function in number theory is defined as follows:

Definition: The classical Smarandache function, $S$, is a function
$S: N \to N$, $N$, *the set of natural numbers such that*
$S(1) = 1$, *and*
$S(n) =$ *The smallest integer such that* $n/S(n)!$

This function has been extensively studied and many interesting properties of it have been discovered [1]. Subsequently many Smarandache type functions have been defined and their interesting properties have been achieved. Ruiz and Perez have discussed some properties of several Smarandache type functions that are involved in many proposed, solved and unsolved problems [2].

An exact formula for counting primes less than or equal to given x in terms of classical Smarandache function has been discovered by L. Seagull [3]. Ruiz and Perez have quoted this result along with a proof while discussing some properties of the classical Smarandache function (Property 2.4) [2].

## 2. A formula for twin prime pairs:

We now proceed to obtain an exact formula for counting twin prime pairs less than or equal to given x in terms of the classical Smarandache function.

We denote by $T_2(x)$ the exact number of twin prime pairs less than or equal to x. Also [m] denotes the integral part of m.

**Theorem:**

$$T_2(x) = -1 + \sum_{1}^{x-2}\left[\frac{S(j)S(j+2)}{(j)(j+2)}\right]$$

where, $S(k)$ denotes the value of classical Smarandache function evaluated at k.

Proof: It is well known that
  (i) $S(p) = p$ iff p is prime $> 4$
  (ii) $S(p) < p$ when p is not prime and $p \neq 4$.
  (iii) $S(4) = 4$.

*In the light of the above properties,*

$$\left[\frac{S(2)S(4)}{(2)(4)}\right] = 1$$

therefore, (2, 4) will be counted as a twin prime pair in the sum given in the above formula. The term " -1 " is added in the formula to eliminate this additional count. Also,

$$\left[\frac{S(j)S(j+2)}{(j)(j+2)}\right] = 1$$

Only when $(j, j+2)$ will be a twin prime pair and in all other cases

$$\left[\frac{S(j)S(j+2)}{(j)(j+2)}\right] = 0.$$

Hence the theorem.

Let us denote by $T_{2n}(x)$ the exact number of prime pairs $(p, p+2n)$, n a positive integer and n > 1.

Corollary:

$$T_{2n}(x) = \sum_{1}^{x-2n}\left[\frac{S(j)S(j+2n)}{(j)(j+2n)}\right]$$

Proof of the corollary: Since n > 1, the illegal appearance of the pair (2, 4) as a prime pair is automatically prohibited, and the proof follows by proceeding on the similar lines.

3. Conclusion: Like formula for counting primes up to given x, [3], one can obtain a similar formula for counting twin prime pairs as well as prime pairs in which the primes are separated by *2n* in terms of the Classical Smarandache Function.

## Acknowledgements

The author is thankful to Dr. M. R. Modak and Dr. S. A. Katre, Bhaskaracharya Pratishthana, Pune, for their keen interest.